\documentclass[conference]{IEEEtran}
\IEEEoverridecommandlockouts
\usepackage{cite}
\usepackage{amsmath,amssymb,amsfonts,amsthm}
\usepackage{algorithmic}
\usepackage{graphicx}
\usepackage{textcomp}
\usepackage{xcolor}
\usepackage{multirow}
\usepackage{caption}
\usepackage{subcaption}
\usepackage{url}
\def\BibTeX{{\rm B\kern-.05em{\sc i\kern-.025em b}\kern-.08em
    T\kern-.1667em\lower.7ex\hbox{E}\kern-.125emX}}
    
\theoremstyle{definition}
\newtheorem{definition}{Definition}
\newtheorem{proposition}{Proposition}
\newtheorem{theorem}{Theorem}
\newcommand{\MCP}{\text{MCP}}

\begin{document}

\title{A mixed complementarity problem approach for steady-state voltage and frequency stability analysis\\
\thanks{This  material  was  based upon work supported by the U.S. Department of Energy, Office of Science, under Contract No. DE-AC02-06CH11357.}
}

\author{\IEEEauthorblockN{Youngdae Kim}
\IEEEauthorblockA{\textit{Mathematics and Computer Science Division} \\
\textit{Argonne National Laboratory}\\
Lemont, IL, USA\\
youngdae@anl.gov}
\and
\IEEEauthorblockN{Kibaek Kim}
\IEEEauthorblockA{\textit{Mathematics and Computer Science} \\
\textit{Argonne National Laboratory}\\
Lemont, IL, USA\\
kimk@anl.gov}
}

\maketitle

\begin{abstract}
We present a mixed complementarity problem (MCP) approach for a steady-state stability analysis of voltage and frequency of electrical grids.
We perform a theoretical analysis providing conditions for the global convergence and local quadratic convergence of our solution procedure, enabling fast computation time. 
Moreover, algebraic equations for power flow, voltage control, and frequency control are compactly and incrementally formulated as a single MCP that subsequently is solved by a highly efficient and robust solution method.
Experimental results over large grids demonstrate that our approach is as fast as the existing Newton method with heuristics while showing  more robust performance.
\end{abstract}


\section{Introduction}
\label{sec:intro}

The growing penetration of renewable energy resources into the electrical grid and the increasing interdependency between the natural gas network and the grid draw attention to the impact of these components on grid stability.
Unlike traditional generators, the intermittency and lack of reactive power generation capability of renewable energy resources can significantly degrade grid stability. Natural gas as the largest source of U.S. electricity generation~\cite{EIA21} also can have a huge impact on stable operations of the grid, as we have seen in  the Texas power outage in the winter storm of early 2021.

One of the key measures for assessing grid stability is a steady-state power flow analysis with voltage regulation and frequency control.
It basically computes a solution to alternating current power flows of the grid formulated as a nonlinear system of equations  that also satisfies several complicated conditions modeling the interactions between grid components to implement those regulations.
An example of such conditions is voltage set point control via reactive power, which regulates voltage magnitude to stay at its set point by controlling its corresponding reactive power.

Computationally, these additional conditions make it difficult to directly apply the the Newton--Raphson (NR) method~\cite{Tinney-Hart67}, which has been the norm for  conventional power flow analysis without regulations.
Heuristics have been developed~\cite{Stott74,Zhao08} that fix some variable values and switch the bus type (e.g., PV-PQ switching) at each iteration of the NR method so that it can still be used to solve  the reduced system of equations.
However, they can potentially cause numerical issues~\cite{Zhao08} leading to divergence.
The lack of theoretical analysis of conditions that guarantee the convergence of these heuristics further complicates the issue to avoid such divergence.

Recent works~\cite{Sundaresh14,Murray15,Tinoco16} introduce complementarity constraints
$0 \le a \perp b \ge 0$ (which denotes $ab=0,a,b \ge 0$) 
to directly incorporate regulation conditions in the model.
They then transform the problem into a computationally more tractable form via equation reformulation using the Fischer--Burmeister function defined by $\phi(a,b)=\sqrt{a^2+b^2}-a-b$ so that $0 \le a \perp b \ge 0 \Leftrightarrow \phi(a,b)=0$.
In particular,~\cite{Murray15} showed global convergence results that lay grounds for a numerically more robust performance than that of the NR method.
The computation time of this approach, however, is usually much slower than that of the NR method with heuristics, as also shown in~\cite{Murray15,Tinoco16}.
The reason is partly that the solution procedure does not exploit complementarity structure since it is obscured by reformulations into equations, resulting in no special treatment possible for it.

In this paper we show that such problems can be cast as a mixed complementarity problem (\MCP{}) with a proper formulation, as will be described in Section~\ref{sec:formulations}.
In particular, we can leverage the convergence theory and solution methods for \MCP{}s to efficiently tackle our problem.
We show sufficient conditions for obtaining global convergence as well as local quadratic convergence that enable fast computation time.
Such fast convergence is made through the solution method for \MCP{}~\cite{Dirkse95,Ferris99-2} that specifically utilizes the complementarity structure.
Similar to the NR method, at each iteration it linearizes the problem at the current iterate; however, instead of solving a system of equations, it solves a linear complementarity problem (LCP), which provides a first-order approximation to the original \MCP{}.
Numerical experiments in Section~\ref{sec:exp} demonstrate that the computational time of our \MCP{} approach is as fast as that of the NR method, while showing a more robust performance.

The rest of the paper is organized as follows.
In Section~\ref{sec:background} we briefly introduce \MCP{} formulations, strong regularity, and their solution method.
Section~\ref{sec:formulations} introduces our \MCP{} formulations for a steady-state power flow analysis with voltage and frequency regulations.
Numerical results are presented in Section~\ref{sec:exp}, and we conclude the paper in Section~\ref{sec:conclusion}.
\section{Background of mixed complementarity problems}
\label{sec:background}

In this section we briefly introduce \MCP{}s, their solution method based on the generalized equation~\cite{Robinson79}, and conditions that guarantee local superlinear or quadratic convergence.

An $\MCP{}(B,F)$ is defined by a vector-valued continuous function $F:\mathbf{R}^n \rightarrow \mathbf{R}^n$ and a box constraint $B=[l,u]=\prod_{i=1}^n[l_i,u_i]$ with $l_i \le u_i$ and $l_i, u_i \in \mathbf{R} \cup \{-\infty,+\infty\}$ for $i=1,\dots,n$.
We note that $F$ is a square system of equations; that is, it has the same number of variables and equations.
We say that $x$ is a solution to the $\MCP{}(B,F)$ if it satisfies one of the following three conditions for each $i=1,\dots,n$:
\begin{equation}
    \begin{aligned}
    x_i = l_i \quad \& \quad F_i(x) \ge 0\\
    l_i \le x_i \le u_i \quad \& \quad F_i(x)=0\\
    x_i=u_i \quad \& \quad F_i(x) \le 0 .
    \end{aligned}
    \label{eq:mcp}
\end{equation}

We use a complementarity notation $l_i \le x_i \le u_i \perp F_i(x)$ to denote~\eqref{eq:mcp}.
Its vector form $l \le x \le u \perp F(x)$ implies that~\eqref{eq:mcp} holds componentwise.\footnote{As a variation, we may use the notation $l_i \le x_i \perp F_i(x) \ge 0$ when $u_i=\infty$. This emphasizes that $F_i(x)$ should be nonnegative at a solution since the third condition of~\eqref{eq:mcp} cannot hold in this case. Similarly, we may denote $x_i \le u_i \perp F_i(x) \le 0$ when $l_i=-\infty$. If $l_i=-\infty$ and $u_i=\infty$, we may omit both bounds and simply denote $x_i \perp F_i(x)$.}

MCPs subsume many different problem classes.
When $B=\mathbf{R}^n$, the $\MCP(B,F)$ becomes a square system of nonlinear equations seeking $x \in \mathbf{R}^n$ that satisfies $F(x)=0$.
The Karush--Kuhn--Tucker conditions of an optimization problem $\min_{l \le x \le u}\; f(x) \; \text{s.t.}\; c(x)=0$ can be formulated as an $\MCP(B,F)$ with $B=[l,u] \times \mathbf{R}^m$ and $F(x,\lambda)=((\nabla f(x)+\nabla c(x)\lambda)^T, (c(x))^T)^T \in \mathbf{R}^{n+m}$, where $c(x),\lambda \in \mathbf{R}^m$.
Many other different problem classes, including generalized Nash equilibrium problems and quasi-variational inequalities, also can be formulated as MCPs. 
We refer to~\cite{Kim18,Kim19-MPC} for  details.

The generalized equation (GE) provides a machinery to solve MCPs~\cite{Dirkse95,Josephy79}.
It is defined by
\begin{equation}
    0 \in F(x) + N_B(x),
    \label{eq:ge}
\end{equation}
where $N_B(x)$ is a normal map to $B$ at $x \in B$ defined by $N_B(x):=\{v \mid \langle v, y-x \rangle \le 0, \forall y \in B\}$ and $N_B(x):=\emptyset$ when $x \notin B$.
One can easily verify that $x$ satisfies~\eqref{eq:ge} if and only if it satisfies~\eqref{eq:mcp}.
A Newton's method~\cite{Dirkse95,Josephy79} finds a solution to the GE by iteratively linearizing and solving the linearized GE, where it becomes an LCP in this case with $F$ being an affine function.
This method is similar to the NR method for a nonlinear system of equations. However, it takes into account complementarity at its linearization, and its complementary pivoting solution procedure effectively exploits the complementarity structure in finding a solution.
We will see its fast and robust computational performance in Section~\ref{sec:exp}.

A key condition to define the domain of attraction for local convergence to a solution $x^*$ of~\eqref{eq:ge} is strong regularity, defined below:
\begin{definition}[Strong regularity~\cite{Robinson80}]
Let $x^*$ be a solution to~\eqref{eq:ge}.
Let $T:=LF_{x^*} + N_B$, where $LF_{x^*}:=F(x^*)+\nabla F(x^*)(x-x^*)$.
We say that $x^*$ is strongly regular if and only if there exist neighborhoods $U$ of the origin and $V$ of $x^*$ such that the restriction to $U$ of $T^{-1}\cap V$ is a single-valued Lipschitzian function from $U$ to $V$.
\label{def:strong-regularity}
\end{definition}
We note that for a system of nonlinear equations $F(x)=0$ strong regularity reduces to the condition that $\nabla F(x^*)$ has a continuous linear inverse map.
This condition is a typical assumption for the local convergence of the NR method.
In Section~\ref{sec:formulations} we present conditions for local convergence of our power flow analysis with voltage and frequency regulations via MCPs.

The following shows an asymptotic convergence rate of the Newton's method~\cite{Dirkse95,Josephy79} for the GE under strong regularity.
\begin{theorem}[Local $Q$-quadratic convergence~\cite{Josephy79}]
For a given $\MCP(B,F)$ with $B \neq \emptyset$, suppose that $x^*$ is a strongly regular solution of it and that $\nabla F$ is locally Lipschitz continuous.
Then for $x^k$ near $x^*$, the sequence generated by solving LCPs is locally Q-quadratically convergent to $x^*$.
\label{thm:local-convergence}
\end{theorem}

\section{Mixed complementarity formulations for a unified power flow analysis}
\label{sec:formulations}

This section presents \MCP{} formulations that encapsulate power flow equations, voltage regulation, and frequency control in an incremental fashion.
Starting with power flow equations in Section~\ref{subsec:power-flow}, we incrementally build \MCP{} formulations in Sections~\ref{subsec:voltage-control}--\ref{subsec:frequency-control}, each of which models a specific regulation and control feature.
At the end of this section, we describe the global and local convergence results.
We note that complementarity constraints for voltage and frequency regulations are known in the literature, but our work is the first showing that these also can be compactly formulated as a single \MCP{}, enabling us to utilize its convergence theory and solution method.

\subsection{Modeling power flow equations}
\label{subsec:power-flow}

\begin{table}[t]
    \caption{Variables fixed/free at each bus $i$ according to its type}
    \label{tbl:fix-vars}
    \centering
    \begin{tabular}{|c|c|c|}
        \hline
         Bus type & Variables fixed & Variables free \\\hline
         Slack & $v_i,\delta_i$ & $p_{g_i},q_{g_i}$\\
         PQ bus & $p_{g_i},q_{g_i}$ & $v_i,\delta_i$\\
         PV bus & $p_{g_i},v_i$ & $\delta_i,q_{g_i}$\\
         \hline
    \end{tabular}
\end{table}

In  conventional power flow analysis, we are interested in finding voltage magnitudes and angles of the grid for a given set of parameters such as load, real power generation, voltage magnitudes at regulated buses, and network coefficients.
To achieve this, we solve the following system of nonlinear equations defined at each bus $i=1,\dots,N$:
\begin{equation}
    \begin{aligned}
        P_i(x) &= p^{\text{inj}}_i - \sum_{k=1}^N v_iv_k(G_{ik}\cos(\delta_{ik})+B_{ik}\sin(\delta_{ik}))\\
        Q_i(x) &= q^{\text{inj}}_i - \sum_{k=1}^N v_iv_k(G_{ik}\sin(\delta_{ik})-B_{ik}\cos(\delta_{ik})),
    \end{aligned}
    \label{eq:pf-def}
\end{equation}
where $p^{\text{inj}}_i := p_{g_i} - P_{d_i}$ and $q^{\text{inj}}_i := q_{g_i} - Q_{d_i}$ are the net real and reactive power injections at bus $i$ with $p_{g_i}$ and $q_{g_i}$ being real and reactive power generated at the bus and $P_{d_i}$ and $Q_{d_i}$ denoting its real and reactive loads.
Voltage magnitude and angle at bus $i$ are denoted by $v_i$ and $\delta_i$, respectively, and the angle difference is defined by $\delta_{ik} := \delta_i - \delta_k$.
We denote by $x:=([p_{g_i},q_{g_i},v_i,\delta_i]_{i=1}^N)$ the encapsulation of all variables.
$G_{ik}$ and $B_{ik}$ are network coefficients computed from a nodal admittance matrix and bus shunt values.
Later in Section~\ref{subsec:transfomer-switched-shunt} they may be changed to variables for voltage regulation.

When solving~\eqref{eq:pf-def}, some components of variable $x$ are fixed depending on the bus type as described in Table~\ref{tbl:fix-vars}.
Since we are interested mainly in voltage magnitudes and angles, we solve for equations that involve these voltage values as variables.
These equations correspond to the real and reactive power flow equations at PQ buses and real power flow equations at PV buses.
Therefore the conventional power flow analysis solves the following square system of nonlinear equations:
\begin{equation}
    \begin{bmatrix}
    P_{\text{PQ}}(x)\\
    Q_{\text{PQ}}(x)\\
    P_{\text{PV}}(x)
    \end{bmatrix} = 0,
    \label{eq:pf-eq}
\end{equation}
where we use the notation $\text{PQ}$ or $\text{PV}$ in the subscript to indicate the set of bus indices belonging to the indicated bus type.
We note that variables in~\eqref{eq:pf-eq} are $\delta_{\text{PQ}}, v_{\text{PQ}}$, and $\delta_{\text{PV}}$.

Using the fact that variables are set to be free variables in~\eqref{eq:pf-eq}, we can formulate~\eqref{eq:pf-eq} as an \MCP{} as follows:
\begin{equation}
    \begin{aligned}
        \delta_{\text{PQ}} &&\perp&& P_{\text{PQ}}(x)\\
        v_{\text{PQ}} &&\perp&& Q_{\text{PQ}}(x)\\
        \delta_{\text{PV}} &&\perp&& P_{\text{PV}}(x)\\.
    \end{aligned}
    \label{eq:mcp-pf}
\end{equation}
From~\eqref{eq:mcp}, one can easily verify that $x$ solves~\eqref{eq:pf-eq} if and only if it solves~\eqref{eq:mcp-pf}.


\subsection{Modeling generator voltage control}
\label{subsec:voltage-control}

Generator voltage control can be compactly formulated by using the following MCP formulation:
\begin{equation}
    \begin{aligned}
        q^{\min}_{g_{\text{PV}}} \le q_{g_{\text{PV}}} \le q^{\max}_{g_{\text{PV}}} &&\perp&& v_{\text{PV}} - v^{\text{sp}}_{\text{PV}},
    \end{aligned}
    \label{eq:mcp-gen-vol}
\end{equation}
where $v^{\text{sp}}_i$ denotes the set point of the voltage magnitude for bus $i \in \text{PV}$.
From~\eqref{eq:mcp}, we see that $v_i=v^{\text{sp}}_i$ when $q^{\min}_{g_i} < q_{g_i} < q^{\max}_{g_i}$, $v_i \ge v^{\text{sp}}_i$ for $q_{g_i}=q^{\min}_{g_i}$, and $v_i \le v^{\text{sp}}_i$ for $q_{g_i}=q^{\max}_{g_i}$.
Therefore~\eqref{eq:mcp-gen-vol} formulates the desired behavior for generator voltage control.
In contrast to the conventional power flow analysis, we note that $v_i, i \in \text{PV}$ in~\eqref{eq:mcp-gen-vol} is not fixed to its set point.

Since both $v_i$ and $q_{g_i}$ are variables for $i \in \text{PV}$, the reactive power flow equations at PV buses should be added in addition to the complementarity constraints~\eqref{eq:mcp-gen-vol}.
Together with~\eqref{eq:mcp-pf}, the following incrementally formulates a power flow analysis with generator voltage control as an \MCP{}:
\begin{equation}
    \begin{aligned}
        v_{\text{PV}} &&\perp&& Q_{\text{PV}}(x)\\
        q^{\min}_{g_{\text{PV}}} \le q_{g_{\text{PV}}} \le q^{\max}_{g_{\text{PV}}} &&\perp&& v_{\text{PV}} - v^{\text{sp}}_{\text{PV}}.
    \end{aligned}
    \label{eq:mcp-gen-vol-full}
\end{equation}



\subsection{Modeling tap-charging transformer and switched shunt voltage control}
\label{subsec:transfomer-switched-shunt}

These additional controllers are used to regulate the bounds on voltage magnitudes as defined below using \MCP{} formulations.
We note that a similar formulation has been given in~\cite{Tinney-Hart67}.
\begin{equation}
    \begin{aligned}
        0 \le u^+_{\text{PV}} &&\perp&& v_{\text{PV}} + v^-_{\text{PV}} - v^{\min}_{\text{PV}} \ge 0\\
        0 \le u^-_{\text{PV}} &&\perp&& v^{\max}_{\text{PV}} - v_{\text{PV}} + v^+_{\text{PV}} \ge 0\\
        0 \le v^+_{\text{PV}} &&\perp&& u_{\text{PV}} - u^{\min}_{\text{PV}} \ge 0\\
        0 \le v^-_{\text{PV}} &&\perp&& u^{\max}_{\text{PV}} - u_{\text{PV}} \ge 0\\
        u_{\text{PV}} &&\perp&& u_{\text{PV}} - u^{\text{sp}}_{\text{PV}} - u^+_{\text{PV}} + u^-_{\text{PV}}
    \end{aligned}
    \label{eq:mcp-tf-ss-vol}
\end{equation}

The first two conditions allow an increase (or a decrease) of a tap ratio $u_{\text{PV}}$ from its set point $u^{\text{sp}}_{\text{PV}}$, when a voltage magnitude reaches its limit.
The third and fourth conditions together with the first two conditions capture a further drift from voltage magnitude limits when the corresponding tap ratio reaches its limit.
For switched shunt devices, we can employ similar \MCP{} formulations.

We note that tap ratios and susceptances in general take discrete values—choosing among finite values—but we instead treat them as continuous variables in~\eqref{eq:mcp-tf-ss-vol}.
This is a common approach as in \cite{Stott74} and still provides useful values for finding a discrete solution~\cite{Tinoco16,Stott74,Chang88}.

In addition to the generator voltage control by~\eqref{eq:mcp-gen-vol-full}, the formulation~\eqref{eq:mcp-tf-ss-vol} can be easily concatenated to~\eqref{eq:mcp-gen-vol-full} to include in the model the effect of tap-changing transformer and switched shunt devices for voltage control as well.
This provides a great flexibility in modeling the impact of various regulation components in power flow analysis via MCPs. 
In Section 4.3, we present a modeling example of incrementally including voltage regulating components to regulate the voltage magnitude.



\subsection{Modeling primary frequency control}
\label{subsec:frequency-control}


Primary frequency control is the governor's reaction to real power imbalance and can be formulated by using an \MCP{} formulation as described in~\eqref{eq:mcp-freq}.
When $p^{\min}_{g_i} < p_{g_i} < p^{\max}_{g_i}$, we have $p_{g_i} = p^{\text{sp}}_{g_i} + \nu_{g_i}\Delta f$ following the linear generation characteristic.
If $p_{g_i}=p^{\max}_{g_i}$, then it stays at its limit; but the frequency could further drop, resulting in $p_{g_i} < p^{\text{sp}}_{g_i} + \nu_{g_i}\Delta f$.
Similarly, we could have $p_{g_i} > p^{\text{sp}}_{g_i} + \nu_{g_i}\Delta f$ when $p_{g_i}=p^{\min}_{g_i}$.
In~\eqref{eq:mcp-freq}, we note that the real power flow equation of the slack bus is added to the formulation for the synchronized generation.

\begin{equation}
    \begin{aligned}
        \Delta f &&\perp&& P_{\text{slack}}(x)\\
        p^{\min}_{g_i} \le p_{g_i} \le p^{\max}_{g_i} &&\perp&& p_{g_i} - (p^{\text{sp}}_{g_i} + \nu_{g_i} \Delta f)\\
        &&&&\forall i \in \text{PV} \cup \text{Slack}
    \end{aligned}
    \label{eq:mcp-freq}
\end{equation}

By combining~\eqref{eq:mcp-freq} with voltage regulation formulations in Sections~\ref{subsec:voltage-control} and~\ref{subsec:transfomer-switched-shunt}, we can perform a power flow analysis with voltage and frequency regulations simultaneously using a single \MCP{} model.


\subsection{Convergence results}

The following proposition shows sufficient conditions for strong regularity of a solution to an \MCP{}.
The conditions can be applied to any \MCP{}s defined in the preceding sections.
By applying Theorem~\ref{thm:local-convergence} with strong regularity, we obtain a local $Q$-quadratic convergence result.
\begin{proposition}[Theorem 3.1 in~\cite{Robinson80}]
Let $x^*$ be a solution to an $\MCP(B,F)$.
Let $\alpha:=\{i \mid x^{\min}_i < x^*_i < x^{\max}_i, F_i(x^*)=0\}, \beta:=\{i \mid x^*_i \in \{x^{\min}_i,x^{\max}_i\}, F_i(x^*)=0\}$, and $\gamma := \{i \mid x^*_i \in \{x^{\min}_i,x^{\max}_i\}, F_i(x^*)\neq 0\}$.
If the submatrix $\nabla F(x^*)_{\alpha,\alpha}$ is nonsingular and the Schur complement $\nabla F(x^*)_{\alpha\cup\beta,\alpha\cup\beta}/\nabla F(x^*)_{\alpha,\alpha}:=\nabla F(x^*)_{\beta,\beta} - \nabla F(x^*)_{\beta,\alpha}\nabla F(x^*)_{\alpha,\alpha}^{-1}\nabla F(x^*)_{\alpha,\beta}$ has positive principal minors (a $P$-matrix), then $x^*$ is a strongly regular solution.
\label{prop:local-conv}
\end{proposition}

We note that when $\beta = \emptyset$ in Proposition~\ref{prop:local-conv}, the conditions correspond to having a continuous linear inverse mapping for the reduced space defined by $\alpha$.

Global convergence results are obtained by assuming a sequence having an accumulation point that converges to a strongly regular solution.
See~\cite[Theorem 5]{Ferris99-2} for details.

\section{Experimental results}
\label{sec:exp}

In this section we present numerical results of our \MCP{} approach and compare its performance with those of the NR method with switching heuristics and the FB-based optimization approach~\cite{Murray15}.
Experiments were performed over large grids included in the MATPOWER package~\cite{Zimmerman11} on a Mac machine having Intel 6-Core i7@2.6 GHz and 32 GB of memory.
The PATH and Ipopt were used for solving \MCP{} and the FB-based problem, respectively.

\subsection{Generator voltage control}
\label{subsec:exp-gen-vol}

In Table~\ref{tbl:exp-gen-vol} we demonstrate the computational performance of each method for voltage control via reactive power.
For \MCP{} we solve the \MCP{} defined by~\eqref{eq:mcp-pf} and~\eqref{eq:mcp-gen-vol-full} in this case.
Our \MCP{} method showed the most robust performance as solving all of the problems, whereas the NR method failed the convergence on 3120sp, and the FB method could not find a solution for ACTIVSg70k.
Also, the \MCP{} method demonstrated  computation time as fast as that of the NR method, while the FB method showed a much slower computation time.

\begin{table*}[!t]
    \centering
    \caption{Generator voltage control}
    \label{tbl:exp-gen}
    \begin{tabular}{|r|r|r|r|r|r|r|r|r|r|}
    \hline
       \multicolumn{1}{|c|}{\multirow{3}{*}{Data}}  &  \multicolumn{3}{c|}{NR method} & \multicolumn{3}{c|}{\MCP{}} & \multicolumn{3}{c|}{FB method}\\\cline{2-10}
       & \multicolumn{1}{c|}{\multirow{2}{*}{Iter}} & \multicolumn{1}{c|}{Time} & \multicolumn{1}{c|}{$\max$}
       & \multicolumn{1}{c|}{\multirow{2}{*}{Iter}} & \multicolumn{1}{c|}{Time} & \multicolumn{1}{c|}{$\max$}
       & \multicolumn{1}{c|}{\multirow{2}{*}{Iter}} & \multicolumn{1}{c|}{Time} & \multicolumn{1}{c|}{$\max$}\\
       & & (secs) & $|v-v^{\text{sp}}|$ & & (secs) & $|v-v^{\text{sp}}|$ & & (secs) & $|v-v^{\text{sp}}|$\\\hline
       1354pegase & 4 & 0.06 & 2.64e-02 & 4 & 0.05 & 2.64e-02 & 22 & 0.82 & 2.64e-02\\
       2869pegase & 6 & 0.18 & 1.82e-02 & 6 & 0.29 & 1.64e-02 & 24 & 2.21 & 1.65e-02\\
       3120sp & f & n/a & n/a & 6 & 0.15 & 6.95e-02 & 925 & 82.15 & 6.95e-02\\
       6468rte & 5 & 0.26 & 2.79e-02 & 4 & 0.28 & 2.78e-02 & 301 & 66.83 & 2.79e-02\\
       9241pegase & 6 & 0.52 & 2.47e-02 & 7 & 0.92 & 2.46e-02 & 74 & 36.95 & 2.47e-02\\
       13659pegase & 5 & 0.63 & 5.66e-02 & 5 & 0.93 & 5.02e-03 & 35 & 23.59 & 5.03e-03\\
       ACTIVSg10k & 5 & 0.36 & 4.07e-05 & 3 & 1.40 & 4.06e-05 & 37 & 10.95 & 1.67e-04\\
       ACTIVSg25k & 4 & 0.80 & 5.83e-04 & 5 & 1.37 & 5.82e-04 & 253 & 263.05 & 7.31e-04\\
       ACTIVSg70k & 5 & 3.12 & 1.04e-03 & 5 & 3.55 & 1.03e-03 & f & n/a & n/a\\\hline
    \end{tabular}
    \label{tbl:exp-gen-vol}
\end{table*}

We found for the 3120sp data that the solutions from the FB method and our \MCP{} method violate some of the voltage magnitude bounds.
Since reactive power cannot control such violations once it reaches its limit, we need another controller such as tap ratios and switched shunt devices, as described in Section~\ref{subsec:transfomer-switched-shunt}, in order to guarantee the feasibility.
This will be addressed in the following section.

\subsection{Voltage control using transformer tap ratios and switched shunt devices}
\label{subsec:exp-tt-ss-vol}

To control the problematic voltage magnitudes of the 3120sp data to be within their bounds, we incorporated~\eqref{eq:mcp-tf-ss-vol} into our \MCP{} problem (so we solved~\eqref{eq:mcp-pf}+\eqref{eq:mcp-gen-vol-full}+\eqref{eq:mcp-tf-ss-vol}) by allowing some buses to change their tap ratios and switched shunt devices.

\begin{figure}[tbhp]
    \centering
    \begin{subfigure}[b]{0.18\textwidth}
    \centering
        \includegraphics[scale=.18]{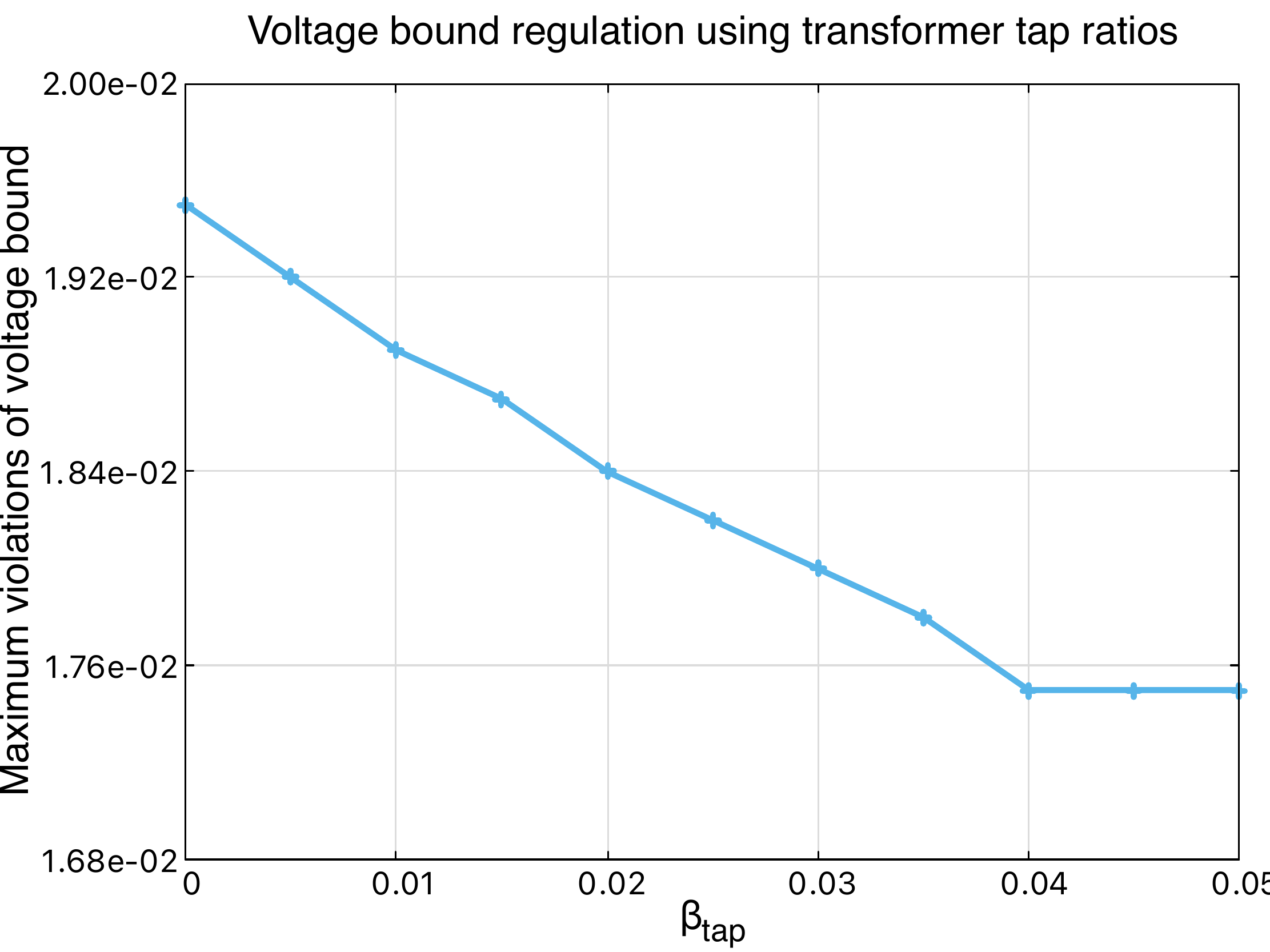}
        \caption{Using transformer tap ratios only}
        \label{subfig:exp-tap}
    \end{subfigure}
    \hspace{8mm}
    \begin{subfigure}[b]{0.18\textwidth}
    \centering
        \includegraphics[scale=.18]{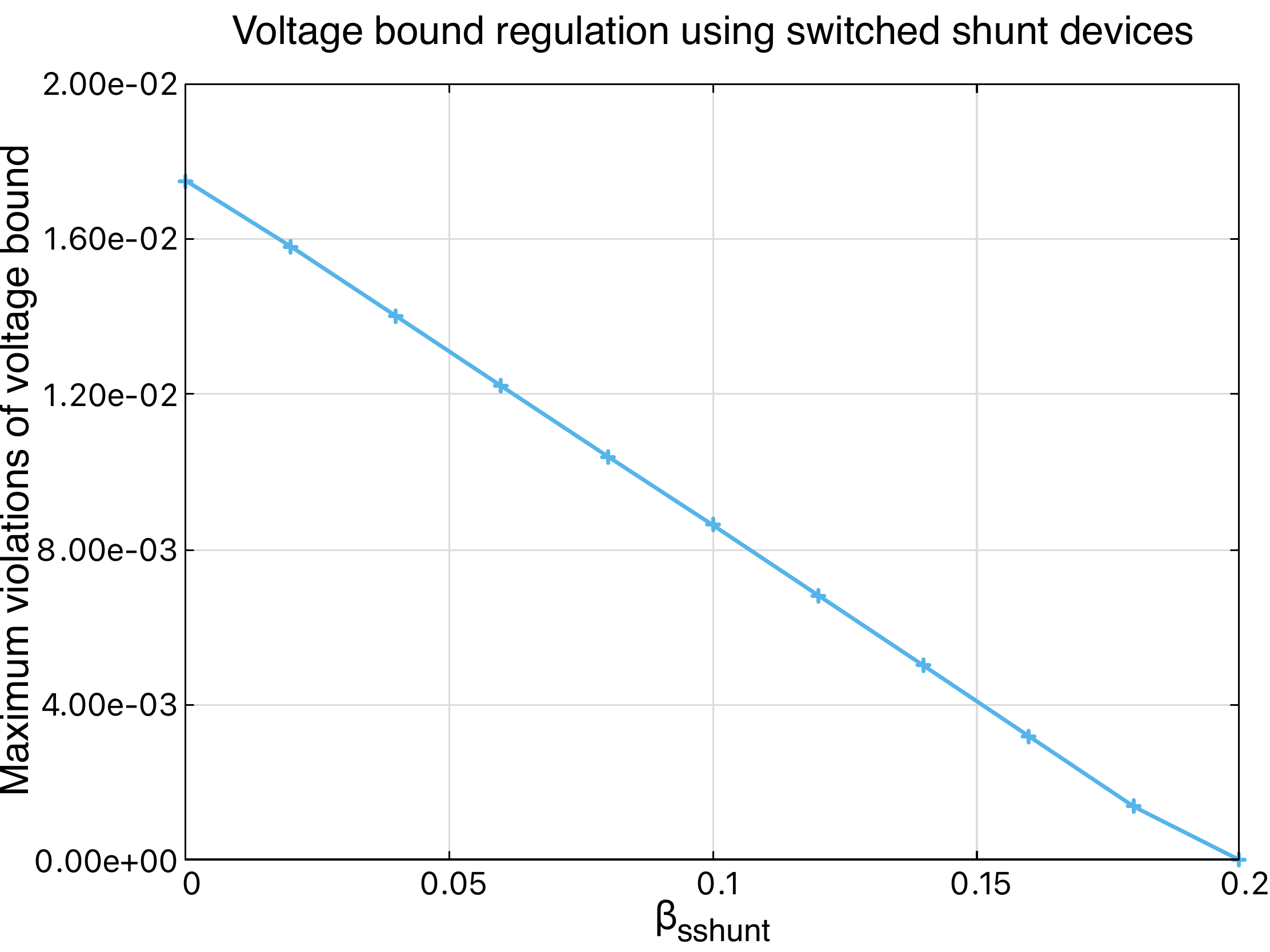}
        \caption{Using both tap ratios and switched shunts}
        \label{subfig:exp-sshunt}
    \end{subfigure}
    \caption{Voltage range regulation using transformers and shunt devices over 3120sp}
    \label{fig:exp-tr-ss-vol}
\end{figure}

Figure~\ref{fig:exp-tr-ss-vol} demonstrates how voltage bound violations were decreasing and eventually removed as we increased the allowed bounds on the tap ratios and switched shunt devices.
These results demonstrate the capability of our method to easily simulate the effect of several voltage controllers in an incremental fashion.

\subsection{Frequency and voltage regulations.}
\label{subsec:exp-freq-vol}

In Table~\ref{tbl:exp-freq}, we report the numerical results from our method for the MCP problem \eqref{eq:mcp-pf}+\eqref{eq:mcp-gen-vol-full}+\eqref{eq:mcp-freq} with frequency and voltage controls on ACTIVSg25k data.
For frequency control we gradually increased real power generation loss by turning off up to 5 generators that had the most real power generation among other generators (i.e., IDs 1557, 2345, 2346, 235, and 2015 in the order of removed).
Those generators were chosen from the solution of generator voltage control we performed in Section~\ref{subsec:exp-gen-vol}, and we warm-started from that solution to solve the \MCP{}.
To the best of our knowledge,  only one recent paper~\cite{sanchez2020integrated}  uses  NR-based heuristics for solving the problem without any convergence guarantee. However, we do not consider NR method because of the lack of convergence guarantee.

In all cases, the \MCP{} was able to quickly find a solution containing frequency changes while controlling the voltage magnitudes to stay at their set points as much as possible and be within their bounds as well.

\begin{table}[tbhp]
    \centering
    \caption{Frequency and voltage regulations over ACTIVSg25k}
    \label{tbl:exp-freq}
    \begin{tabular}{|c|c|c|c|}
    \hline
       Generation loss  & Frequency & $\max |v-v^{\text{sp}}|$ & Time (secs)\\\hline
       1,299 MW & 59.92 Hz & 3.41e-03 & 1.27\\
       2,597 MW & 59.84 Hz & 9.50e-03 & 1.26\\
       3,895 MW & 59.76 Hz & 1.96e-02 & 1.28\\
       5,185 MW & 59.68 Hz & 2.80e-02 & 1.56\\
       6,455 MW & 59.62 Hz & 2.90e-02 & 1.77\\\hline
    \end{tabular}
    \label{tbl:exp-freq}
\end{table}

\section{Conclusion}
\label{sec:conclusion}

We presented a mixed complementarity problem approach for a steady-state power flow analysis with voltage regulation and frequency control. The key advantage of our approach is leveraging the existing algorithms with the theoretical support for the global convergence with quadratic local convergence rate, which guarantees numerical stability and fast computation.
We demonstrated the greater computational performance of our approach, as compared with the existing NR method and FB method, by using large MATPOWER test instances.

\bibliographystyle{IEEEtran}
\bibliography{pes-bib}

 \vspace{3mm}
 {\footnotesize
 \noindent\fbox{\parbox{0.47\textwidth}{
 The submitted manuscript has been created by UChicago Argonne, LLC, Operator of Argonne National Laboratory (``Argonne''). Argonne, a U.S. Department of Energy Office of Science laboratory, is operated under Contract No. DE-AC02-06CH11357. The U.S. Government retains for itself, and others acting on its behalf, a paid-up nonexclusive, irrevocable worldwide license in said article to reproduce, prepare derivative works, distribute copies to the public, and perform publicly and display publicly, by or on behalf of the Government. The Department of Energy will provide public access to these results of federally sponsored research in accordance with the DOE Public Access Plan (http://energy.gov/downloads/doe-public-access-plan).}
 }
 }

\end{document}